\documentclass[12pt,draft]{article}

\usepackage[cp1251]{inputenc}
\usepackage[T2A]{fontenc}
\usepackage[english,ukrainian]{babel}
\usepackage{amsmath,amsfonts,amssymb}
\usepackage{geometry}

\begin{document}

\begin{center}

\large

\textbf{H\"ORMANDER SPACES ON MANIFOLDS, \\ AND THEIR APPLICATION TO ELLIPTIC BOUNDARY-VALUE PROBLEMS}

\medskip

\textbf{T.M. Kasirenko, A.A. Murach, I.S. Chepurukhina}

\medskip

\normalsize

Institute of Mathematics, National Academy of Sciences of Ukraine, Kyiv

\medskip

E-mail: kasirenko@imath.kiev.ua, murach@imath.kiev.ua, Chepurukhina@gmail.com

\bigskip\bigskip

\large

\textbf{ПРОСТОРИ ХЕРМАНДЕРА НА МНОГОВИДАХ ТА\\ ЇХ ЗАСТОСУВАННЯ ДО ЕЛІПТИЧНИХ КРАЙОВИХ ЗАДАЧ}\footnote{Публікація містить результати досліджень, проведених за грантом Президента України за конкурсним проектом Ф75/29007 Державного фонду фундаментальних досліджень.}

\medskip

\textbf{Т.М. Касіренко, О.О. Мурач, І.С. Чепурухіна}

\normalsize

\end{center}

\medskip

\noindent We introduce an extended Sobolev scale on a smooth compact manifold with boundary. The scale is formed by inner-product H\"ormander spaces for which an RO-varying radial function serves as a regularity index. These spaces do not depend on a choice of local charts on the manifold. The scale consists of all Hilbert spaces that are interpolation ones for pairs of inner-product Sobolev spaces, is obtained by the interpolation with a function parameter of these pairs, and is closed with respect to this interpolation. As an application of the scale introduced, we give a theorem on the Fredholm property of a general elliptic boundary-value problem on appropriate H\"ormander spaces and find sufficient conditions under which its generalized solutions belong to the space of $p\geq0$ times continuously differential functions.

\medskip


\bigskip

\noindent \textbf{Вступ.} У сучасному математичному аналізі важливу роль відіграють простори розподілів, для яких показником регулярності служить не число (як у класичних просторах Соболєва), а досить загальний функціональний параметр, залежний від частотних змінних (див., наприклад, [1\,--\,5]). Більше ніж півстоліття тому Л. Хермандер [\ref{Hermander65}] увів і дослідив широкі класи таких просторів та навів їх застосування до диференціальних рівнянь, заданих у евклідових областях. Втім, довгий час простори Хермандера не знаходили широкого застосування у теорії багатовимірних крайових задач, що було пов'язано з браком зручних аналітичних методів для роботи з цими просторами і відсутністю коректного їх означення на многовидах. В останній час ситуація істотно змінилася завдяки роботам В.А. Михайлеця, О.О. Мурача та їх учнів (див. монографію [\ref{MikhailetsMurach14}] і наведену там літературу). Ними виділено класи гільбертових просторів Хермандера, які отримуються інтерполяцією з функціональним параметром пар соболєвських просторів і допускають коректне означення на многовидах (незалежне від вибору локальних карт). Для таких класів вдалося побудувати теорію розв'язності загальних еліптичних крайових задач.

Серед цих класів найширшим є сім'я усіх гільбертових просторів, інтерполяційних для пар гільбертових просторів Соболєва. Її уведено і досліджено в [6\,--\,8] для евклідових областей і замкнених гладких многовидів та названо розширеною соболєвською шкалою. Вона замкнена відносно інтерполяції гільбертових просторів з функціональним параметром.

Мета цієї роботи --- увести розширену соболєвську шкалу на довільному гладкому компактному многовиді з краєм і дослідити її властивості, зокрема, інтерполяційні. Окрім того, ми наведемо деякі застосування цієї шкали до загальних еліптичних крайових задач.

\textbf{1. Простори Хермандера на многовидах.} Нехай $M$ --- компактний орієнтовний нескінченно гладкий многовид вимірності $n\geq1$ з краєм $\partial M\neq\varnothing$. Уведемо клас гільбертових функціональних просторів на $M^\circ:=M\setminus \partial M$, узявши за основу простори Хермандера  $B_{p,\mu}(\mathbb{R}^{n})$ [\ref{Hermander65}, п. 2.2; \ref{Hermander86}, п. 10.1], де число $p=2$, а функція $\mu$ частотного аргументу $\xi\in\mathbb{R}^{n}$ набирає вигляду  $\mu(\xi)=\varphi(\langle\xi\rangle)$; тут $\langle\xi\rangle:=(1+|\xi|^{2})^{1/2}$, а $\varphi\in\mathrm{RO}$.

За означенням, множина $\mathrm{RO}$ складається з усіх вимірних за Борелем функцій $\varphi:\nobreak[1,\infty)\rightarrow(0,\infty)$, для яких існують числа $b>1$ і $c\geq1$ такі, що $c^{-1}\leq\varphi(\lambda t)/\varphi(t)\leq c$ для будь-яких $t\geq1$ і $\lambda\in[1,b]$ (числа $b$ і $c$ залежать від $\varphi$). Клас $\mathrm{RO}$ введений В.Г.~Авакумовичем [\ref{Avakumovic36}], допускає простий опис і досить повно вивчений (див., наприклад, [\ref{Seneta85}, приложение]). Функцію $\varphi\in\mathrm{RO}$ називають RO-змінною на нескінченності.

Як відомо [\ref{Seneta85}, с. 88], для кожної функції $\varphi\in\mathrm{RO}$ існують дійсні числа $s_{0}\leq s_{1}$ та додатні числа $c_{0}$ і $c_{1}$ такі, що
\begin{equation}\label{9f1b}
c_{0}\lambda^{s_{0}}\leq\frac{\varphi(\lambda t)}{\varphi (t)}\leq
c_{1}\lambda^{s_{1}} \quad\mbox{для всіх}\quad t\geq1,\;\;\lambda\geq1.
\end{equation}
Поклавши тут $t:=1$, бачимо, що ця функція є міжстепеневою. Її зв'язок зі степеневими функціями характеризують числа $\sigma_{0}(\varphi)$ і $\sigma_{1}(\varphi)$, перше з яких є супремумом усіх дійсних $s_{0}$ таких, що виконується ліва частина нерівності \eqref{9f1b}, а друге є інфімумом усіх дійсних $s_{1}$ таких, що виконується права частина нерівності \eqref{9f1b}. Числа
$\sigma_{0}(\varphi)$ і $\sigma_{1}(\varphi)$ називаються відповідно нижнім і верхнім індексами Матушевської [\ref{Matuszewska64}] функції $\varphi\in\mathrm{RO}$. Зокрема, якщо вона правильно змінна на нескінченності [\ref{Seneta85}, с. 9], то ці індекси дорівнюють її порядку.

Нехай $\varphi\in\mathrm{RO}$. Нагадаємо означення вказаного простору Хермандера $B_{2,\varphi(\langle\cdot\rangle)}(\mathbb{R}^{n})$, який будемо позначати через $H^{\varphi}(\mathbb{R}^{n})$. Ми розглядаємо комплекснозначні функції і розподіли та комплексні функціональні простори. За означенням, лінійний простір
$H^{\varphi}(\mathbb{R}^{n})$ складається з усіх повільно зростаючих розподілів $w$ на $\mathbb{R}^{n}$ таких, що їх перетворення Фур'є $\widehat{w}$ є локально інтегровним за Лебегом на $\mathbb{R}^{n}$ і задовольняє умову
\begin{equation*}
\|w\|_{\varphi,\mathbb{R}^{n}}^{2}:=\int_{\mathbb{R}^{n}}
\varphi^2(\langle\xi\rangle)\,|\widehat{w}(\xi)|^2\,d\xi<\infty.
\end{equation*}
Простір $H^{\varphi}(\mathbb{R}^{n})$ гільбертів і сепарабельний відносно норми $\|\cdot\|_{\varphi,\mathbb{R}^{n}}$.

У випадку степеневої функції $\varphi(t)\equiv t^{s}$ він стає гільбертовим простором Соболєва $H^{(s)}(\mathbb{R}^{n})$ порядку $s\in\mathbb{R}$. Узагалі виконуються щільні неперервні вкладення $H^{(s_{1})}(\mathbb{R}^{n})\hookrightarrow H^{\varphi}(\mathbb{R}^{n})\hookrightarrow H^{(s_{0})}(\mathbb{R}^{n})$ для довільних дійсних чисел $s_{0}$ і $s_{1}$, які задовольняють умову \eqref{9f1b}, зокрема, для довільних $s_{0}<\sigma_{0}(\varphi)$ і $s_{1}>\sigma_{1}(\varphi)$. Клас просторів $\{H^{\varphi}(\mathbb{R}^{n}):\varphi\in\mathrm{RO}\}$ досліджено в [\ref{MikhailetsMurach13UMJ3}; \ref{MikhailetsMurach14}, п. 2.4.2] і названо розширеною соболєвською шкалою.

Для відкритої непорожньої множини $\Omega\subset\mathbb{R}^{n}$ простір
$H^{\varphi}(\Omega)$ складається, за означенням, із звужень на $\Omega$ усіх розподілів $w\in H^{\varphi}(\mathbb{R}^{n})$. Цей простір гільбертів і сепарабельний відносно норми
$$
\|u\|_{\varphi,\Omega}:=\inf\bigl\{\|w\|_{\varphi,\mathbb{R}^{n}}:
w\in H^{\varphi}(\mathbb{R}^n),\;\,u=w\;\,\mbox{в}\;\,\Omega\bigr\},
$$
де $u\in H^{\varphi}(\Omega)$. Він є ізотропним випадком просторів, досліджених Л.Р. Волєвичем і Б.П. Панеяхом [\ref{VolevichPaneah65}, \S~3]. Нас окремо буде цікавити випадок, коли $\Omega$ є півпростором $\mathbb{R}^n_+:=\{(x',x_n): x'\in\mathbb{R}^{n-1},\;x_n>0\}$.

Означимо тепер гільбертів простір $H^{\varphi}(M^\circ)$ за допомогою локальних карт і розбиття одиниці на $M$ та норм у просторах $H^{\varphi}(\mathbb{R}^{n})$ і $H^{\varphi}(\mathbb{R}^n_+)$.
Із $C^{\infty}$-структури на компактному многовиді $M$ виберемо який-небудь його скінченний атлас. Не обмежуючи загальності, вважаємо, що останній складається з $\varkappa$ локальних карт
$\pi_j:\overline{\mathbb{R}^n_+}\leftrightarrow M_j$, де $j=1,\ldots,\varkappa$, і $\varkappa_{0}$ локальних карт
$\pi_j:\mathbb{R}^n\leftrightarrow M_j$, де $j=\varkappa+1,\ldots,\varkappa+\varkappa_{0}$. Тут відкриті (у топології на $M$) множини $M_j$, де $j=1,\ldots,\varkappa+\varkappa_{0}$, утворюють покриття многовиду $M$ таке, що $\nobreak{M_j\cap\partial M\neq\varnothing}$ тоді і лише тоді, коли $j\leq\varkappa$. Звісно, $\overline{\mathbb{R}^n_+}:=\{(x',x_n): x'\in\mathbb{R}^{n-1},\;x_n\geq\nobreak0\}$. Окрім того, виберемо розбиття одиниці на $M$, утворене деякими функціями $\chi_j\in C^{\infty}(M)$, де $j=1,\ldots,\varkappa+\varkappa_{0}$, які задовольняють умову $\mathrm{supp}\,\chi_j\subset M_j$.

Нехай, як і раніше, $\varphi\in\mathrm{RO}$. За означенням, простір
$H^{\varphi}(M^\circ)$, є поповненням лінійного многовиду $C^{\infty}(M)$ за нормою
$$
\|u\|_{\varphi,M^\circ}=\biggl(\,\sum_{j=1}^{\varkappa}
\|(\chi_j u)\circ\pi_j\|^2_{\varphi,\mathbb{R}^n_+}+
\sum_{j=\varkappa+1}^{\varkappa+\varkappa_{0}}
\|(\chi_j u)\circ\pi_j\|^2_{\varphi,\mathbb{R}^n}\biggr)^{1/2}.
$$
Тут, звісно, $(\chi_j u)\circ\pi_j$ позначає функцію $(\chi_j u)(\pi_j(x))$ аргументу $x\in\overline{\mathbb{R}^n_+}$ або $x\in\mathbb{R}^n$. Цей простір гільбертів і сепарабельний відносно вказаної норми.

\textbf{Теорема 1.} \it Гільбертів простір $H^{\varphi}(M^\circ)$, де $\varphi\in\mathrm{RO}$, не залежить з точністю до еквівалентності норм від вибору атласу многовиду $M$ і відповідного розбиття одиниці на $M$. \rm

Клас просторів $\{H^{\varphi}(M^\circ):\varphi\in\mathrm{RO}\}$ називаємо розширеною соболєвською шкалою на $M^\circ$. Якщо $\varphi(t)\equiv t^{s}$ для деякого $s\in\mathbb{R}$, то $H^{\varphi}(M^\circ)$ стає гільбертовим простором Соболєва порядку $s$, який позначаємо через $H^{(s)}(M^\circ)$.
У випадку, коли функція $\varphi$ правильно змінна на нескінченності за Й.~Караматою, простір $H^{\varphi}(M^\circ)$, уведено і досліджено в [\ref{MikhailetsMurach06UMJ3}, п. 3]. Якщо $M^\circ$~--- евклідова область $\Omega$ (обмежена з межею класу $C^{\infty}$), то $H^{\varphi}(M^\circ)=H^{\varphi}(\Omega)$ з точністю до еквівалентності норм. У цьому випадку розширена соболєвська шкала досліджена в [\ref{MikhailetsMurach15}] (навіть для областей з ліпшіцевою межею).

Кожний простір $H^{\varphi}(M^\circ)$ неперервно вкладається у лінійний топологічний простір усіх продовжуваних розподілів на $M^\circ$. Тому ці простори можна порівнювати.

\textbf{Теорема 2.} \it Нехай $\varphi_{0},\varphi_{1}\in\mathrm{RO}$. Вкладення $H^{\varphi_1}(M^\circ)\subset H^{\varphi_0}(M^\circ)$ виконується тоді і тільки тоді, коли функція $\varphi_0/\varphi_1$ обмежена в околі нескінченності. У цьому випадку вкладення неперервне. Воно компактне тоді і тільки тоді, коли $\varphi_0(t)/\varphi_1(t)\to0$ при $t\to\infty$. \rm

З теореми 2 випливає, що для довільних дійсних чисел $s_{0}$ і $s_{1}$, які задовольняють \eqref{9f1b}, виконуються неперервні вкладення $H^{(s_{1})}(M^\circ)\hookrightarrow H^{\varphi}(M^\circ)\hookrightarrow H^{(s_{0})}(M^\circ)$. Ці вкладення компактні, якщо $s_{0}<\sigma_{0}(\varphi)$ і $s_{1}>\sigma_{1}(\varphi)$.

Обговоримо зв'язок простору $H^{\varphi}(M^\circ)$ з його аналогом на $\partial M$; тут $n\geq2$. Зауважимо, що $\partial M$ --- замкнений нескінченно гладкий многовид вимірності $n-1$. Простір $H^{\eta}(\partial M)$, де $\eta\in\mathrm{RO}$, уведено і досліджено в [\ref{MikhailetsMurach09Dop3}] (див. також [\ref{MikhailetsMurach14}, п. 2.4.2]). Він є поповненням лінійного многовиду $C^{\infty}(\partial M)$ за нормою
$$
\|h\|_{\eta,\partial M}=\biggl(\,\sum_{j=1}^{\varkappa}
\|(\chi_{j}h)(\pi_j(\cdot,0))\|^2_{\eta,\mathbb{R}^{n-1}}
\biggr)^{1/2}.
$$
Простір $H^{\eta}(\partial M)$ гільбертів та сепарабельний відносно цієї норми і з точністю до еквівалентності норм не залежить від вибору локальних карт $\pi_j(\cdot,0):\mathbb{R}^{n-1}\leftrightarrow M_j\cap\partial M$, які покривають многовид $\partial M$, та відповідного розбиття одиниці на $\partial M$ [\ref{MikhailetsMurach09Dop3}, с. 32].

\textbf{Теорема 3.} \it Нехай $\varphi\in\mathrm{RO}$ і $\sigma_{0}(\varphi)>1/2$. Тоді відображення $R:u\mapsto u\!\upharpoonright\!\partial M$, де $u\in\nobreak C^\infty(M)$, продовжується єдиним чином (за неперервністю) до обмеженого лінійного оператора сліду $R:H^{\varphi}(M^\circ)\to H^{\varphi\rho^{-1/2}}(\partial M)$. Цей оператор сюр'єктивний і має обмежений лінійний правий обернений оператор $S:H^{\varphi\rho^{-1/2}}(\partial M)\to H^{\varphi}(M^\circ)$ такий, що відображення $S$ не залежить від $\varphi$. \rm

Тут і далі використано функціональний параметр $\varrho(t):=t$ аргументу $t\geq1$. Отже, $\varphi\rho^{-1/2}$ позначає функцію $\varphi(t)t^{-1/2}$ аргументу $t$.

З теореми 3 випливає, що простір $H^{\eta}(\partial M)$, де  $\eta\in\mathrm{RO}$ і $\sigma_{0}(\eta)>0$, складається зі слідів на $\partial M$ усіх розподілів $u\in H^{\eta\rho^{1/2}}(M^\circ)$, а норма у цьому просторі еквівалентна нормі
\begin{equation*}
\inf\bigl\{\|u\|_{\eta\rho^{1/2},M^\circ}:
u\in H^{\eta\rho^{1/2}}(M^\circ),\; Ru=h\bigr\},
\end{equation*}
де $h\in H^\eta(\partial M)$.

\textbf{2. Інтерполяційні властивості просторів Хермандера.} Розширена соболєвська шкала на $M^\circ$ складається (з точністю до еквівалентності норм) з усіх гільбертових просторів, інтерполяційних для пар соболєвських просторів $H^{(s_0)}(M^\circ)$ і $H^{(s_1)}(M^\circ)$, де $s_0<s_1$. Нагадаємо, що гільбертів простір $X$ називають інтерполяційним для пари $[H_0,H_1]$ гільбертових просторів, другий з яких неперервно вкладений у перший, якщо задовольняються такі дві властивості: а) виконуються неперервні вкладення $H_1\hookrightarrow X\hookrightarrow H_0$, б) як тільки який-небудь лінійний оператор $T$ є обмеженим на $H_0$ і на $H_{1}$, то він є також обмеженим на $X$. Зазначена інтерполяційна властивість цієї шкали випливає з такого результату:

\textbf{Теорема 4.} \it Нехай $s_{0},s_{1}\in\mathbb{R}$ і $s_{0}<s_{1}$. Гільбертів простір $X$ є інтерполяційним для пари соболєвських просторів $H^{(s_{0})}(M^\circ)$ і $H^{(s_{1})}(M^\circ)$ тоді і тільки тоді, коли $X=\nobreak H^{\varphi}(M^\circ)$ з точністю до еквівалентності норм для деякого параметра $\varphi\in\mathrm{RO}$, який задовольняє умову \eqref{9f1b}. \rm

Зауважимо, що умову \eqref{9f1b} можна переформулювати за допомогою індексів Матушевської. А саме, вона еквівалентна такій парі умов: i) $s_{0}\leq\sigma_{0}(\varphi)$, та, окрім того,
$s_{0}<\sigma_{0}(\varphi)$, якщо в означенні $\sigma_{0}(\varphi)$ супремум не досягається, ii) $\sigma_{1}(\varphi)\leq s_{1}$ та, окрім того, $\sigma_{1}(\varphi)<s_{1}$, якщо в означенні $\sigma_{1}(\varphi)$ інфімум не досягається.

Для застосувань просторів $H^{\varphi}(M^\circ)$ важливо, що вони отримується інтерполяцією з функціональним параметром деяких пар гільбертових соболєвських просторів на $M^\circ$. У~цьому зв'язку нагадаємо означення методу інтерполяції з функціональним параметром пар гільбертових просторів, запропонованого Ч. Фояшом і Ж.-Л. Ліонсом [\ref{FoiasLions61}]. При його викладі спираємося в основному на монографію [\ref{MikhailetsMurach14}, п. 1.1].

Нам достатньо обмежитися випадком регулярної пари $H:=[H_0,H_1]$ сепарабельних гільбертових просторів. Її регулярність означає, що $H_1$ неперервно і щільно вкладено в $H_0$. Для неї існує самоспряжений додатно визначений оператор $J$ у гільбертовому просторі $H_{0}$ з областю визначення $H_{1}$, який встановлює ізометричний ізоморфізм між гільбертовими просторами $H_{1}$ і $H_{0}$. Цей оператор визначається за парою $H$ однозначно.

Нехай вимірна за Борелем функція $\psi:\nobreak(0,\infty)\rightarrow(0,\infty)$ обмежена на кожному відрізку $[a,b]$, де $0<a<b<\infty$, і відокремлена від нуля на кожній множині $[r,\infty)$, де $r>0$. Множину усіх таких функцій позначимо через $\mathcal{B}$. У гільбертовому  просторі $H_{0}$ за допомогою спектральної теореми  означений (взагалі кажучи, необмежений) оператор $\psi(J)$ як борелева функція $\psi$ від самоспряженого оператора $J$. Позначимо через $[H_{0},H_{1}]_\psi$ або коротко через $H_{\psi}$
область визначення оператора $\psi(J)$, наділену нормою $\|u\|_{H_\psi}:=\|\psi(J)u\|_{H_0}$, де $u\in H_\psi$. Простір $H_\psi$ гільбертів і сепарабельний, до того ж виконується неперервне і щільне вкладення $H_\psi\hookrightarrow H_0$.

Функцію $\psi$ називають інтерполяційним параметром, якщо для довільних регулярних пар $H:=[H_0,H_1]$ і $G:=[G_0,G_1]$ гільбертових просторів та для довільного лінійного відображення $T$, заданого на $H_0$, виконується така властивість: якщо при кожному $j\in\{0,1\}$ звуження відображення $T$ на простір $H_{j}$ є обмеженим оператором $T:H_{j}\to G_{j}$, то і звуження відображення $T$ на
простір $H_\psi$ є обмеженим оператором $T:H_\psi\to G_{\psi}$. Тоді кажуть, що простір $H_\psi$ отримано інтерполяцією з функціональним параметром $\psi$ пари $H$. Функція $\psi\in\mathcal{B}$ є інтерполяційним параметром тоді і тільки тоді, коли існує угнута функція $\psi_{0}:(r,\infty)\to(0,\infty)$ така, що обидві функції $\psi/\psi_{0}$ і $\psi_{0}/\psi$ обмежені на $(r,\infty)$, де $r\gg1$.

\textbf{Теорема 5.} \it Нехай $\varphi\in\mathrm{RO}$, а дійсні числа $s_{0}<s_{1}$ задовольняють умову \eqref{9f1b}. Означимо інтерполяційний параметр $\psi\in\mathcal{B}$ за формулами $\psi(\tau):=\tau^{-s_{0}/(s_{1}-s_{0})}\,\varphi(\tau^{1/(s_{1}-s_{0})})$, якщо $\tau\geq1$, і $\psi(\tau):=\varphi(1)$, якщо $0<\tau<1$. Тоді
$$
H^{\varphi}(M^\circ)=
\bigl[H^{(s_{0})}(M^\circ),H^{(s_{1})}(M^\circ)\bigr]_{\psi}
\quad\mbox{з еквівалентністю норм.}
$$\rm

Зауважимо, що будь-які числа $s_{0}<\sigma_{0}(\varphi)$ і $s_{1}>\sigma_{1}(\varphi)$ задовольняють умову \eqref{9f1b} стосовно $\varphi\in\mathrm{RO}$.

Розширена соболєвська шкала на $M^\circ$ замкнена відносно розглянутого методу інтерполяції з функціональним параметром.

\textbf{Теорема 6.} \it Нехай $\varphi_{0},\varphi_{1}\in\mathrm{RO}$ і $\psi\in\mathcal{B}$. Припустимо, що функція $\varphi_{0}/\varphi_{1}$ обмежена в околі нескінченності, а функція $\psi$ є інтерполяційним параметром. Тоді
$$
\bigl[H^{\varphi_{0}}(M^\circ),H^{\varphi_{1}}(M^\circ)\bigr]_{\psi}=
H^{\varphi}(M^\circ)\quad\mbox{з еквівалентністю норм},
$$
де функція $\varphi(t):=\varphi_{0}(t)\,\psi(\varphi_{1}(t)/\varphi_{0}(t))$ аргументу $t\geq1$ належить до класу $\mathrm{RO}$. \rm

У випадку, коли $M^\circ$~--- евклідова область, теореми 4\,--\,6 доведено в [\ref{MikhailetsMurach15}].

\textbf{3. Застосування.} Розглянемо на $M^\circ$ еліптичну крайову задачу
\begin{equation}\label{9f7}
Au=f\quad\mbox{на}\;\,M^\circ,\qquad
B_{j}u=g_{j}\quad\mbox{на}\;\,\partial M,
\quad j=1,...,q
\end{equation}
(див., наприклад, [\ref{Agranovich97}, п. 1.2]). Тут $A$ --- лінійний диференціальний оператор на $M$ довільного парного порядку $2q\geq2$, а кожне $B_{j}$~--- крайовий лінійний диференціальний оператор на $\partial M$ довільного порядку $m_{j}\geq0$. Усі коефіцієнти цих операторів належать до класів $C^{\infty}(M)$ і $C^{\infty}(\partial M)$ відповідно. Покладемо $B:=(B_{1},\ldots,B_{q})$ і $m:=\mathrm{max}\{m_{1},\ldots,m_{q}\}$.

\textbf{Теорема 7.} \it Припустимо, що $\varphi\in\mathrm{RO}$ і
$\sigma_0(\varphi)>m+1/2$. Тоді відображення $u\mapsto\nobreak(Au,B_{1}u,,\ldots,B_{q}u)$, де $u\in C^{\infty}(M)$, продовжується єдиним чином (за неперервністю) до обмеженого лінійного оператора на парі гільбертових просторів
\begin{equation*}
H^{\varphi}(M^\circ)\quad\mbox{і}\quad H^{\varphi\varrho^{-2q}}(M^\circ)\oplus\bigoplus_{j=1}^{q}
H^{\varphi\varrho^{-m_j-1/2}}(\partial M).
\end{equation*}
Цей оператор нетерів. Його ядро лежить в $C^{\infty}(M)$ і разом з індексом не залежить від $\varphi$.\rm

Ця теорема виводиться з соболєвського випадку (коли $\varphi(t)\equiv t^{s}$) за допомогою інтерполяційної теореми 5 та її аналогу для просторів Хермандера на $\partial M$.

Нехай $U$ --- відкрита (у топології на $M$) підмножина многовиду $M$, а $V:=U\cap\partial M\neq\varnothing$. Позначимо через $H^{\varphi}_{\mathrm{loc}}(U)$, де $\varphi\in\mathrm{RO}$, лінійний простір усіх продовжуваних розподілів $u$ на $M^\circ$ таких, що $\chi u\in H^{\varphi}(M^\circ)$ для довільної функції $\chi\in C^{\infty}(M)$, яка задовольняє умову $\mathrm{supp}\,\chi\subset U$. Аналогічно, позначимо через $H^{\eta}_{\mathrm{loc}}(V)$, де $\eta\in\mathrm{RO}$, лінійний простір усіх розподілів $h$ на $\partial M$ таких, що $\chi h\in H^{\eta}(\partial M)$ для довільної функції $\chi\in C^{\infty}(\partial M)$, яка задовольняє умову $\mathrm{supp}\,\chi\subset V$.

\textbf{Теорема 8.} \it Припустимо, що функція $u\in\bigcup_{s>m+1/2}H^{(s)}(M)$ є розв'язком еліптичної крайової задачі \eqref{9f7}, праві частини якої задовольняють умову
$$
(f,g_{1},\ldots,g_{q})\in H^{\varphi\varrho^{-2q}}_{\mathrm{loc}}(U)\times\prod_{j=1}^{q}
H^{\varphi\varrho^{-m_j-1/2}}_{\mathrm{loc}}(V)
$$
для деякого $\varphi\in\mathrm{RO}$. Тоді $u\in H^{\varphi}_{\mathrm{loc}}(U)$. \rm

З цієї теореми і версії теореми вкладення Хермандера [\ref{Hermander65}, с. 59] для простору $H^{\varphi}(M)$ випливає

\textbf{Теорема 9.} \it Нехай ціле число $p\geq0$. Припустимо, що функція $u$ задовольняє умову теореми $8$ для деякого параметра $\varphi\in\mathrm{RO}$ такого, що $\sigma_{0}(\varphi)>m+1/2$ і
$$
\int_1^{\infty}t^{2p+n-1}\varphi^{-2}(t)dt<\infty.
$$
Тоді $u\in C^{p}(U)$.\rm


\bigskip

\noindent REFERENCES

\begin{enumerate}

\item\label{Hermander65}
H\"ormander, L. (1963). Linear partial differential operators. Berlin: Springer.

\item\label{Hermander86}
H\"ormander, L. (1983). The analysis of linear partial differential operators, vol. II, Differential operators with constant coefficients. Berlin: Springer.

\item\label{Jacob010205}
Jacob, N. (2001, 2002, 2005). Pseudodifferential operators and Markov processes (in 3 volumes). London: Imperial College Press.

\item\label{NicolaRodino10}
Nicola, F. \& Rodino, L. (2010). Global Pseudodifferential Calculus on Euclidean spaces. Basel: Birkh\"aser.

\item\label{MikhailetsMurach14}
Mikhailets, V.A. \& Murach, A.A. (2014). H\"ormander spaces, interpolation, and elliptic problems. Berlin, Boston: De Gruyter.

\item\label{MikhailetsMurach09Dop3}
Mikhailets, V.A. \& Murach, A.A. (2009). Elliptic operators on a closed compact manifold. Dopov. Nac. akad. nauk. Ukr., No. 3, pp. 13-19 (in Russian).

\item\label{MikhailetsMurach13UMJ3}
Mikhailets, V.A. \& Murach, A.A. (2013). Extended Sobolev scale and elliptic operators. Ukr. Math. J., 65, No. 3, pp. 435-447.

\item\label{MikhailetsMurach15}
Mikhailets, V.A. \& Murach, A.A. (2015). Interpolation Hilbert spaces between Sobolev spaces. Results Math, 67, No. 1, pp. 135-152.

\item\label{Avakumovic36}
Avakumovi\'c, V.G. (1936). O jednom O-inverznom stavu. Rad
Jugoslovenske Akad. Znatn. Umjetnosti, 254, pp. 167-186.

\item\label{Seneta85}
Seneta, E. (1976). Regularly varying functions. Berlin: Springer.

\item\label{Matuszewska64}
Matuszewska, W. (1964). On a generalization of regularly increasing functions. Studia Math., 24, pp. 271-279.

\item\label{VolevichPaneah65}
Volevich, L.R. \& Paneah B.P. (1965). Certain spaces of generalized functions and embedding  theorems. Russian Math. Surveys, 20, No. 1, pp. 1-73.

\item\label{MikhailetsMurach06UMJ3}
Mikhailets, V.A. \& Murach, A.A. (2006). Refined scales of spaces, and elliptic boundary value problems. II. Ukr. Math. J., 58, No. 3, pp. 398-417.

\item\label{FoiasLions61}
Foia\c{s}, C. \& Lions, J.-L. (1961). Sur certains th\'eor\`emes d'interpolation. Acta Sci. Math. (Szeged), 22, No. 3-4, pp. 269-282.

\item\label{Agranovich97}
Agranovich, M.S. (1997). Elliptic boundary problems. Encycl. Math. Sci. Vol. 79. Partial differential equations, IX. Berlin: Springer.

\end{enumerate}

\end{document}